\documentclass[a4paper,11pt]{article}
\usepackage[latin1]{inputenc}
\usepackage[T1]{fontenc}
\usepackage[francais]{babel}
\usepackage{amssymb}
\usepackage{amsfonts, amsmath}
\newtheorem{thm}{Theorem}
\newtheorem{prop}{Proposition}

\newtheorem{rem}{Remark}

\title{On Poincar\'{e} lemma or Volterra theorem about differential forms and cohomology groups}
\author{\textbf{A. Lesfari}
\\\emph{Department of Mathematics}
\\\emph{Faculty of Sciences}
\\\emph{University of Choua\"{i}b Doukkali}
\\\emph{B.P. 20, El Jadida, Morocco}.
\\\emph{lesfariahmed@yahoo.fr, lesfari.a@ucd.ac.ma}}
\date{}

\begin{document}
\maketitle \textbf{Abstract.} The Poincar\'{e} lemma (or Volterra
theorem) is of utmost importance both in theory and in practice.
It tells us every differential form which is closed, is locally
exact. In other words, on a contractible manifold all closed forms
are exact. The aim of this paper is to present some direct proofs
of this lemma and explore some of its numerous consequences. Some
connections with Cech-De Rham-Dolbeault cohomologies,
$\overline{\partial}$-Poincar\'{e} lemma or Dolbeault-Grothendieck
lemma are given.

\section{Introduction}

Let $\omega$ be a $k$-differential form ($k>0$) on an open subset
$U\subset\mathbb{R}^n$. We assume that $\omega$ is
$\mathcal{C}^1$. Recall that $\omega$ is closed (or is a cocycle)
if  $d\omega=0$. Similarly, $\omega$ is  exact  (or cohomological
at 0) if there is a $\mathcal{C}^1$, $(k-1)$-differential form
$\lambda$ on $U$ such that : $\omega=d\lambda$. Notice that there
are no exact $0$-forms because there are no $-1$-forms. It is well
known that every exact differential form is closed and that the
reciprocal is false in general. Poincar\'{e} lemma ensures that it
is true on a star-shaped open subset of $\mathbb{R}^n$ whose
definition is as follows : let $[a,b]=\{\lambda a+(1-\lambda)b,
\lambda\in[0,1]\}$, $(a,b)\in\mathbb{R}^n\times\mathbb{R}^n$. We
say that an open subset $U\subset\mathbb{R}^n$ is star-shaped if
there exists an $a$ in $U$, such that for all $b$ in $U$,
$[a,b]\subset U$. In other words, if the line segment from $b$ to
$a$ is in $U$.

We will study Poincar\'{e} lemma which is considered fundamental
both in theory and in practice. Already at the university level
and as evidenced by several scientific books, this is a key result
for the study of many problems in mathematics, physics and the
theorem itself has applications in areas ranging from
electrodynamics to differential and integral calculus on
varieties. At a higher level it intervenes for example when
studying the cohomology of De Rham varieties [7] to name only this
striking example. Let's note for information that the so currently
called Poincar\'{e} lemma is due to Volterra. Indeed, the
Poincar\'{e} lemma is really Volterra theorem; the work of
Volterra is contained in several notes published in the Rendiconti
of the Accademia dei Lincei [10] (see also [6]). In this work we
use as everyone the name Poincar\'{e} lemma instead of Volterra
theorem and we leave the question to be clarified by historians.

The lemma to be demonstrated is a typical example of a local
result, so it suffices to prove it in local coordinates, for
example in an arbitrarily small open of $\mathbb{R}^n$, but it
must be admitted that the technical details are much more
complicated to demonstrate when we move from the case of $1$-forms
to $k$-differential forms. All known proofs of this lemma, using
these coordinates are often too computerized and rarely
illuminating.

We offer some proofs of the Poincar\'{e} lemma in this paper
because the proofs represent widely different views of the
subject. The first proof is given in a very classical setting and
represents the classical point of view; it is both elementary and
constructive but the problem is that it is a bit long and
technical. Then, our problem is to give a quick proof of this
lemma although requiring certain knowledges pushed in differential
geometry; the proofs uses very modern machinery and represent a
more modern point of view. Both of these points of view have merit
and so we demonstrate them both. The paper is divided in some
sections and subsections, each of them devoted to various and
complementary aspects of the problem concerning, in particular,
connections with Cech-De Rham-Dolbeault cohomologies,
$\overline{\partial}$-Poincar\'{e} lemma or Dolbeault-Grothendieck
lemma.

\section{The Poincar\'{e} lemma}

We give a little information (which will be needed in the second
proof) about one-parameter group of diffeomorphisms, differential
operators, Lie derivative, inner product and Cartan's formula.
This discussion is brief, but should be enough to define notation.
Let $M$ be a differentiable manifold of dimension $m$. Let $TM$ be
the tangent bundle to $M$, i.e., the union of spaces tangent to
$M$ at all points $x$, $TM=\bigcup_{x\in M}T_{x}M$. This bundle
has a natural structure of differentiable variety of dimension
$2m$ and it allows us to convey immutably to the manifolds the
whole theory of ordinary differential equations. A vector field
(we also say section of the tangent bundle) on $M$ is an
application, denoted $X$, which at every point $x\in M$ associates
a tangent vector $X_{x}\in T_{x}M$. In other words, $X :
M\longrightarrow TM$, is an application such that if $\pi
:TM\longrightarrow M$, is the natural projection, we have
$\pi\circ X=id_{M}$. Let $\left( x_{1},...,x_{m}\right)$ be a
local coordinate system in a neighborhood $U\subset M$. In this
system the vector field $X$ is written in the form
$$X=\sum_{k=1}^{m}f_{k}(x)\frac{\partial}{\partial x_{k}},
\text{ }x\in U,$$ where the functions
$f_{1},\ldots,f_{m}:U\longrightarrow \mathbb{R}$, are the
components of $X$ with respect to $\left(x_{1},...,x_{m}\right)$.
A vector field $X$ is differentiable if its components $f_{k}(x)$
are differentiable functions. Given a point $x\in M$, we write
$g_{t}^{X}(x)$ (or simply $g_{t}(x)$) the position of $x$ after a
displacement of a duration $t\in \mathbb{R}$. There is thus an
application $g_{t}^{X} : M\longrightarrow M$, $t\in \mathbb{R}$,
which is a diffeomorphism (a one-to-one differentiable mapping
with a differentiable inverse), by virtue of the theory of
differential equations. The vector field $X$ generates a
one-parameter group of diffeomorphisms $g_{t}^{X}$ on $M$, i.e., a
differentiable application ($\mathcal{C}^{\infty }$) : $M\times
\mathbb{R}\longrightarrow M$, satisfying a group law :

$i)$ $\forall t\in \mathbb{R},\text{ }g_{t}^{X} : M\longrightarrow
M$ is a diffeomorphism.

$ii)$ $\forall t,s\in \mathbb{R},\text{
}g_{t+s}^{X}=g_{t}^{X}\circ g_{s}^{X}$.\\
The condition $ii)$ means that the mapping $t\longmapsto
g_{t}^{X}$, is a homomorphism of the additive group $\mathbb{R}$
into the group of diffeomorphisms of $M$. It implies that
$$g_{-t}^{X}=\left(g_{t}^{X}\right)^{-1},$$ because
$g_{0}^{X}=id_{M}$ is the identical transformation that leaves
every point invariant. The one-parameter group of diffeomorphisms
$g_{t}^{X}$ on $M$, which we have just described is called a
flow\index{flow} and it admits the vector field $X$ for velocity
fields
$$\frac{d}{dt}g_{t}^{X}(x)=X\left(g_{t}^{X}(x)\right),$$
with the initial condition : $g_{0}^{X}(x)=x$. Obviously
$$\left.\frac{d}{dt}g_{t}^{X}(x)\right|_{t=0}=X(x).$$
Hence by these formulas $g_{t}^{X}(x)$ is the curve on the
manifold that passes through $x$ and such that the tangent at each
point is the vector $X\left(g_{t}^{X}(x)\right)$. The vector field
$X$ generates a unique group of diffeomorphisms of $M$. With every
vector field $X$ we associate the first-order differential
operator $L_{X}$. This is the differentiation of functions in the
direction of the vector field $X$. We have
$$L_{X}:\mathcal{C}^{\infty}(M)
\longrightarrow \mathcal{C}^{\infty}(M),\text{ }F\longmapsto
L_{X}F,$$ where
$$L_{X}F(x)=\left. \frac{d}{dt}F\left( g_{t}^{X}(x)\right) \right| _{t=0},
\text{ }x\in M.$$ $\mathcal{C}^{\infty }\left( M\right)$ being the
set of functions $F:M\longrightarrow \mathbb{R}$, of class
$\mathcal{C}^{\infty }$. The operator $L_{X}$ is linear :
$L_{X}\left(\alpha_{1}F_{1}+\alpha_{2}F_{2}\right)=\alpha_{1}L_{X}F_{1}+\alpha_{2}L_{X}F_{2}$,
where $\alpha_{1},\alpha_{2}\in \mathbb{R}$, and satisfies the
Leibniz formula :
$$L_{X}\left(F_{1}F_{2}\right)=F_{1}L_{X}F_{2}+F_{2}L_{X}F_{1}.$$
Since $L_{X}F(x)$ only depends on the values of $F$ in the
neighborhood of $x$, we can apply the operator $L_{X}$ to the
functions defined only in the neighborhood of a point, without the
need to extend them to the full variety $M$. Let
$\left(x_{1},...,x_{m}\right)$ be a local coordinate system on
$M$. In this system the vector field $X$ has components $f_{1},
\ldots, f_{m}$ and the flow $g_{t}^{X}$ is defined by a system of
differential equations. Therefore, the derivative of the function
$F=F\left(x_{1}, ...,x_{m}\right)$ in the direction of $X$ is
written
$$L_{X}F=f_{1}\frac{\partial F}{\partial
x_{1}}+\cdots+f_{m}\frac{\partial F}{\partial x_{m}}.$$ In other
words, in the coordinates $\left(x_{1},...,x_{m}\right)$ the
operator $L_{X}$ is written
$\displaystyle{L_{X}=f_{1}\frac{\partial}{\partial
x_{1}}+\cdots+f_{m}\frac{\partial}{\partial x_{m}}}$, this is the
general form of the first order linear differential operator.

Let $M$ and $N$ be two differentiable manifolds of dimension $m$
and $n$ respectively, $U\subset M$, $V\subset N$ two open subsets.
For any differentiable application $g : U\longrightarrow V$, and
any $k$-differential form in $V$,
$$\omega=\sum_{1\leq i_1,...,i_k\leq n}f_{i_1,...,i_k}dx_{i_1}\wedge...\wedge
dx_{i_k},$$ we define a $k$-differential form in $U$ (called the
pull-back by $g$ or inverse image or the transpose of $\omega$ by
$g$) by setting
$$g^*\omega =\sum_{1\leq i_1,...,i_k\leq n}\left(f_{i_1,...,i_k}\circ
g\right)dg_{i_1}\wedge...\wedge dg_{i_k},$$ where
$$dg_{i_l}=\sum_{j=1}^{m}\frac{\partial
g_{i_l}}{\partial y_j}dy_j,$$ are $1$-forms in $U$. Note that
$g^*$ is a linear operator from the space of $k$-forms on $N$ to
the space of $k$-forms on $U$ (the asterisk indicates that $g^*$
operates in the opposite direction of $g$). Let $X$ be a vector
field on a differentiable manifold $M$. We recalled above that the
vector field $X$ generates a unique group of diffeomorphisms
$g_t^X$ (that we also note $g_t$) on $M$, solution of the
differential equation
$$\frac{d}{dt}g_t^X(p)=X(g_t^X(p)),\quad p\in M,$$
with the initial condition $g_0^X(p)=p$. Let $\omega$ be a
$k$-differential form. The Lie derivative of $\omega$ with respect
to $X$ is the $k$-form differential defined by
$$
L_X\omega=\left.\frac{d}{dt}g_t^*\omega\right|_{t=0}=\lim_{t\rightarrow
0}\frac{g_t^*(\omega(g_t(p)))-\omega(p)}{t}.$$ In general, for
$t\neq0$, we have
\begin{equation}\label{eqn:euler}
\frac{d}{dt}g_t^*\omega=\left.\frac{d}{ds}g_{t+s}^*\omega\right|_{s=0}=
g_t^*\left.\frac{d}{ds}g_s^*\omega\right|_{s=0}=g_t^*(L_X\omega).
\end{equation}
It is easily verified that for the $k$-differential form
$\omega(g_t (p))$ at the point $g_t(p)$, the expression
$g_t^*\omega(g_t(p))$ is indeed a $k$-differential form in $p$.
For all $t\in \mathbb{R}$, the application
$g_t:\mathbb{R}\longrightarrow\mathbb{R}$ being a diffeomorphism
then $dg_t$ and $dg_{-t}$ are applications,
\begin{eqnarray}
dg_t&:&T_pM\longrightarrow T_{g_t(p)}M,\nonumber\\
dg_{-t}&:&T_{g_t(p)}M\longrightarrow T_pM.\nonumber
\end{eqnarray}
The Lie derivative of a vector field $Y$ in the direction $X$ is
defined by
$$
L_XY=\left.\frac{d}{dt}g_{-t}Y\right|_{t=0}=\lim_{t\rightarrow
0}\frac{g_{-t}(Y(g_t(p)))-Y(p)}{t}.$$ In general, for $t\neq0$, we
have
$$
\frac{d}{dt}g_{-t}Y=\left.\frac{d}{ds}g_{-t-s}Y\right|_{s=0}
=g_{-t}\left.\frac{d}{ds}g_{-s}Y\right|_{s=0}=g_{-t}(L_Y).$$

An interesting operation on differential forms is the inner
product that is defined as follows : the inner product of a
$k$-differential form $\omega$ by a vector field $X$ on a
differentiable manifold $M$ is a $(k-1)$-differential form,
denoted $i_X\omega$, defined by
$$(i_X\omega)(X_1,...,X_{k-1})=\omega(X,X_1,...,X_{k-1}),$$
where
$X_1,...,X_{k-1}$ are vector fields. If
$$X=\sum_{j=1}^mX_j(x)\frac{\partial}{\partial
x_j},$$
is the local expression of the vector field on the variety
$M$ of dimension $m$ and
$$\omega=\sum_{i_1<i_2<...<i_k}f_{i_1...i_k}(x)dx_{i_1}\wedge...\wedge dx_{i_k},$$
then, the expression of the inner product in local coordinates is
given by
$$i_{\frac{\partial}{\partial
x_j}}\omega=\frac{\partial}{\partial(dx_j)}\omega,$$ where we put
$dx_j$ in first position in $\omega$. Moreover, the well-known
Cartan homotopy formula is fundamental to the Lie derivative and
can be used as a definition :
\begin{equation}\label{eqn:euler}
L_X\omega=d(i_X\omega)+i_X(d\omega).
\end{equation}
For a differential form $\omega$, we have
$$i_XL_X\omega=L_Xi_X\omega.$$

\begin{thm}
On a star-shaped open subset $U$ of $\mathbb{R}^n$, any closed
differential form is exact.
\end{thm}
\emph{\textbf{Proof 1}}: a) (For a $1$-differential form). Let
$U\subset\mathbb{R}^n$, be a star-shaped open subset relative to
one of its points that we note $a$. Let
$$\omega=\sum_{i=1}^nf_idx_i$$,
be a closed $1$-differential form in $U$, and we will show that it
is exact. In other words, that there is a $0$-differential form in
$U$, i.e., a $\mathcal{C}^1$ application
$h:U\longrightarrow\mathbb{R}$ such that : $\omega=dh$ or what
amounts to the same, such that :
$$f_i=\frac{\partial h}{\partial x_i}.$$ (Indeed, if $\omega$ is exact, then by definition
$$\omega=\sum_{i=1}^n f_idx_i=dh=\sum_{i=1}^n\frac{\partial h}{\partial x_i}dx_i,$$
hence the result). We will see that for the regularity of this
$1$-form, $h$ is $\mathcal{C}^2$. By hypothesis, $U$ is a
star-shaped open subset and according to the fundamental theorem
of differential and integral calculus as well as that of
derivation of composite functions, we have
\begin{eqnarray}
h(x)-h(a)&=&\int_0^1\frac{d}{dt}h(a+t(x-a))dt,\nonumber\\
&=&\int_0^1\sum_{j=1}^n\frac{\partial h}{\partial
x_j}(a+t(x-a))).(x_j-a_j)dt.\nonumber
\end{eqnarray}
A function $h$ satisfying the relation
$f_i=\displaystyle{\frac{\partial h}{\partial x_i}}$, is only
defined to an additive constant. We can therefore put
$$h(x)=\int_0^1\sum_{j=1}^n\frac{\partial
h}{\partial x_j}(a+t(x-a))).(x_j-a_j)dt.$$  The application $h$ is
well defined on $U$ because for all $t\in[0,1]$, $a+t(x-a)\in U$
($U$ is a star-shaped open subset). Notice that
$\displaystyle{\frac{\partial h}{\partial x_i}}$ exist; it is a
function defined by an integral. We have
\begin{eqnarray}
\frac{\partial h}{\partial
x_i}(x)&=&\int_0^1\frac{\partial}{\partial
x_i}\left(\sum_{j=1}^nf_j(a+t(x-a)).(x_j-a_j)\right)dt,\nonumber\\
&=&\int_0^1\sum_{j=1}^n\left(\frac{\partial f_j}{\partial
x_i}(a+t(x-a)).t.(x_j-a_j)+f_j(a+t(x-a)).\delta_{ij}\right)dt,\nonumber
\end{eqnarray}
where $\delta_{ij}=1$ if $i=j$ and $0$ if $i\neq j$. Since $w$ is
closed, then
$$
d\omega=\sum_{i=1}^3df_i\wedge dx_i=\sum_{\underset{i<j}{1\leq
i,j\leq 3}}^3\left(\frac{\partial f_j}{\partial
x_i}-\frac{\partial f_i}{\partial x_j}\right)dx_i\wedge dx_j.
$$
As a result, $\omega$ is closed (i.e., $d\omega=0$) if and only if
$$\frac{\partial f_i}{\partial x_j}=\frac{\partial f_j}{\partial
x_i},$$ because $dx_i\wedge dx_j\neq0$, $i\neq j$. If $i=j$, the
relations in question are trivial and for $i>j$, it is obviously
enough to swap the indices $i$ and $j$. Then,
\begin{eqnarray}
\frac{\partial h}{\partial
x_i}(x)&=&\int_0^1\left(t\sum_{j=1}^n\frac{\partial f_i}{\partial
x_j}(a+t(x-a))(x_j-a_j)+f_i(a+t(x-a))\right)dt,\nonumber\\
&=&\int_0^1\left(t\frac{d}{dt}f_i(a+t(x-a))+f_i(a+t(x-a))\right)dt,\nonumber\\
&=&\int_0^1\frac{d}{dt}\left(tf_i(a+t(x-a))\right)dt,\nonumber\\
&=&f_i(x),\nonumber
\end{eqnarray}
where $1\leq i\leq n$ and $x\in U$. The $0$-differential forms in
$U$, i.e., continuous applications $U\longrightarrow\mathbb{R}$
and since $\displaystyle{\omega=\sum_{i=1}^nf_idx_i}$, the
applications $f_i$ are $\mathcal{C}^1$ and we deduce from the
relations above that $h$ is $\mathcal{C}^2$ on $U$.

b) (For a $k$-differential form, $k\geq2$). Without restricting
the generality, we suppose that $U$ is a star-shaped open subset
with respect to $0$. Let $\psi$ be the application from the set of
$l$-differential forms on $U$ to the set of $(l-1)$-differential
forms on $U$ defined by
$$\psi(\lambda)=\sum_{1\leq i_1,...,i_k\leq
n}\sum_{j=1}^l(-1)^{j-1}\left(\int_0^1t^{l-1}g_{i_1,...,i_l}(t.)dt\right)
\pi_{i_j}dx_{i_1}\wedge...\wedge\widehat{dx_{i_j}}\wedge...\wedge
dx_{i_l},$$ where
$$\lambda=\sum_{1\leq i_1,...,i_l\leq n}g_{i_1,...,i_l}dx_{i_1}\wedge...\wedge
dx_{i_l},$$
$$g_{i_1,...,i_l}(t.): U\longrightarrow\mathbb{R},\quad  x\longmapsto g_{i_1,...,i_l}(tx),
t\in[0,1],$$ the application $\pi_{i_j}$ is the projection on the
$i_j^{\mbox{th}}$ coordinate and $\widehat{dx_{i_j}}$ denotes the
term omitted. Let
$$\omega=\sum_{1\leq i_1,...,i_k\leq n}f_{i_1,...,i_k}dx_{i_1}\wedge...\wedge
dx_{i_k},$$ be a $\mathcal{C}^1$, $k$-differential form in $U$.
The idea of the proof is to show that
$$\omega=\psi(d\omega)+d(\psi(\omega)),$$
because in this case since by hypothesis $\omega$ is closed, then
$\omega=d(\psi(\omega))$ (by construction, $\psi(0)=0$) and
therefore, the form $\omega$ is exact (note that $\psi(\omega)$ is
$\mathcal{C}^1$ if $\omega$ is). We have
$$d\omega=\sum_{1\leq i_1,...,i_k\leq n}\frac{\partial f_{i_1,...,i_k}}{\partial x_i}\wedge dx_i\wedge dx_{i_1}\wedge...\wedge
dx_{i_k},$$ and
\begin{eqnarray}
&&\psi(d\omega)\nonumber\\
&&=\sum_{i=1}^n\sum_{1\leq i_1,...,i_k\leq n}\left[
\left(\int_0^1t^k\frac{\partial f_{i_1,...,i_k}}{\partial
x_i}(t.)dt\right)p_idx_{i_1}\wedge...\wedge
dx_{i_k}\right.\nonumber\\
&&\quad+\sum_{j=1}^k(-1)^j\left.\left(\int_0^1t^k\frac{\partial
f_{i_1,...,i_k}}{\partial x_i}(t.)dt\right)p_{i_j}dx_i\wedge
dx_{i_1}\wedge...\wedge\widehat{dx_{i_j}}\wedge...\wedge
dx_{i_k}\right],\nonumber
\end{eqnarray}
or
\begin{eqnarray}
&&\psi(d\omega)\nonumber\\
&&=\sum_{1\leq i_1,...,i_k\leq n}\left[\sum_{i=1}^n
\left(\int_0^1t^k\frac{\partial f_{i_1,...,i_k}}{\partial
x_i}(t.)dt\right)p_idx_{i_1}\wedge...\wedge
dx_{i_k}\right]\nonumber\\
&&\quad-\sum_{1\leq i_1,...,i_k\leq
n}\left[\sum_{i=1}^n\sum_{j=1}^k(-1)^{j-1}\left(\int_0^1t^k\frac{\partial
f_{i_1,...,i_k}}{\partial x_i}(t.)dt\right)\right.\nonumber\\
&&\left.\qquad\qquad\qquad\qquad\qquad\qquad\qquad\qquad
p_{i_j}dx_i\wedge
dx_{i_1}\wedge...\wedge\widehat{dx_{i_j}}\wedge...\wedge
dx_{i_k}\right].\nonumber
\end{eqnarray}
Similarly, we have
\begin{eqnarray}
&&d(\psi(\omega))\nonumber\\
&&=\sum_{1\leq i_1,...,i_k\leq
n}\sum_{j=1}^k(-1)^{j-1}d\left[\left(\int_0^1t^{k-1}
f_{i_1,...,i_k}(t.)dt\right)p_{i_j}\right]\nonumber\\
&&\qquad\qquad\qquad\qquad\qquad\qquad\qquad\qquad\qquad\quad\wedge
dx_{i_1}\wedge...\wedge\widehat{dx_{i_j}}\wedge...\wedge
dx_{i_k},\nonumber
\end{eqnarray}
or
\begin{eqnarray}
&&d(\psi(\omega))\nonumber\\ &&=\sum_{1\leq i_1,...,i_k\leq
n}\sum_{j=1}^k(-1)^{j-1}\sum_{i=1}^n\left[\left(\int_0^1t^{k-1}
\frac{\partial f_{i_1,...,i_k}}{\partial
x_i}(t.)tdt\right)p_{i_j}\right.\nonumber\\
&&\quad\qquad+\left.\left(\int_0^1t^{k-1}
f_{i_1,...,i_k}(t.)dt\right)\delta_{i_j,i}\right]dx_i\wedge
dx_{i_1}\wedge...\wedge\widehat{dx_{i_j}}\wedge...\wedge
dx_{i_k},\nonumber
\end{eqnarray}
or
\begin{eqnarray}
&&d(\psi(\omega))\nonumber\\&&=\sum_{1\leq i_1,...,i_k\leq
n}\sum_{i=1}^n\sum_{j=1}^k(-1)^{j-1}\left(\int_0^1t^k
\frac{\partial f_{i_1,...,i_k}}{\partial x_i}(t.)dt\right)\nonumber\\
&&\qquad\qquad\qquad\qquad\qquad\qquad\qquad\qquad
p_{i_j}dx_i\wedge
dx_{i_1}\wedge...\wedge\widehat{dx_{i_j}}\wedge...\wedge dx_{i_k}\nonumber\\
&&\qquad+\sum_{1\leq i_1,...,i_k\leq
n}\sum_{j=1}^k(-1)^{j-1}\left(\int_0^1t^{k-1}
f_{i_1,...,i_k}(t.)dt\right)\nonumber\\
&&\qquad\qquad\qquad\qquad\qquad\qquad\qquad\qquad dx_{i_j}\wedge
dx_{i_1}\wedge...\wedge\widehat{dx_{i_j}}\wedge...\wedge
dx_{i_k},\nonumber
\end{eqnarray}
and finally,
\begin{eqnarray}
&&d(\psi(\omega))\nonumber\\
&&= \sum_{1\leq i_1,...,i_k\leq
n}\sum_{i=1}^n\sum_{j=1}^k(-1)^{j-1}\left(\int_0^1t^k
\frac{\partial f_{i_1,...,i_k}}{\partial x_i}(t.)dt\right)\nonumber\\
&&\qquad\qquad\qquad\qquad\qquad\qquad\qquad\qquad
p_{i_j}dx_i\wedge
dx_{i_1}\wedge...\wedge\widehat{dx_{i_j}}\wedge...\wedge dx_{i_k}\nonumber\\
&&\qquad+k\sum_{1\leq i_1,...,i_k\leq n}\left(\int_0^1t^{k-1}
f_{i_1,...,i_k}(t.)dt\right)dx_{i_1}\wedge...\wedge
dx_{i_k}.\nonumber
\end{eqnarray}
Hence,
\begin{eqnarray}
&&\psi(d\omega)+d(\psi(\omega))\nonumber\\
&&= \sum_{1\leq i_1,...,i_k\leq n}\sum_{i=1}^n\left(\int_0^1t^k
\frac{\partial f_{i_1,...,i_k}}{\partial x_i}(t.)dt\right)p_i
dx_{i_1}\wedge...\wedge dx_{i_k}\nonumber\\
&&\qquad+\sum_{1\leq i_1,...,i_k\leq n}\left(\int_0^1kt^{k-1}
f_{i_1,...,i_k}(t.)dt\right)dx_{i_1}\wedge...\wedge dx_{i_k},\nonumber\\
&&=\sum_{1\leq i_1,...,i_k\leq
n}\left(\int_0^1\frac{d}{dt}\left(t^k
f_{i_1,...,i_k}(t.)\right)dt\right)dx_{i_1}\wedge...\wedge dx_{i_k},\nonumber\\
&&=\sum_{1\leq i_1,...,i_k\leq n}
f_{i_1,...,i_k}dx_{i_1}\wedge...\wedge dx_{i_k},\nonumber\\
&&=\omega,\nonumber
\end{eqnarray}
and proof 1 ends.

\emph{\textbf{Proof 2}}: Using the notions and properties
mentioned at the beginning of this section, we give a quick proof
of the lemma in question. Indeed, consider the differential
equation
$$\dot{x}=X(x)=\frac{x}{t},$$ as well as its solution
$$g_t(x_0)=x_0t.$$
The latter is defined in the neighborhood of the point $x_0$,
depending on how $\mathcal{C}^\infty$ of the initial condition and
is a parameter group of diffeomorphisms. We have,
$$g_0(x_0)=0,\quad
g_1(x_0)=x_0,\quad g_0^*\omega=0,\quad g_1^*\omega=\omega.$$
Hence,
\begin{eqnarray}
\omega&=&g_1^*\omega-g_0^*\omega,\nonumber\\
&=&\int_0^1\frac{d}{dt}g_t^*\omega dt,\nonumber\\
&=&\int_0^1g_t^*(L_X\omega)dt\quad(\mbox{by (1)}),\nonumber\\
&=&\int_0^1g_t^*(di_X\omega)dt,\quad(\mbox{according to (2) and the fact that }d\omega=0)\nonumber\\
&=&\int_0^1dg_t^*i_X\omega dt,\quad(\mbox{because }
df^*\omega=f^*d\omega).\nonumber
\end{eqnarray}
We can therefore find a differential form $\lambda$ such that :
$\omega=d\lambda$, where $$\lambda=\int_0^1g_t^*i_X\omega dt,$$
which completes the proof 2 and ends the two proofs of the lemma.
$\square$

\begin{rem}
Any exact differential form is closed. It is well known that the
converse is false in general and depends on the open $U$ on which
the differential form is $\mathcal{C}^1$. For example, if
$U=\mathbb{R}^2\backslash\{(0,0)\}$ then the differential form
$$\omega=-\frac{x_2}{x_1^2+x_2^2}dx_1+\frac{x_1}{x_1^2+x_2^2}dx_2,$$ is closed but is not exact.
Indeed, this form is obviously closed. To show that it is not
exact, we use the fact that in general if $\omega$ is an exact
$1$-differential form on an open subset and $\gamma$ a closed path
in this $\mathcal{C}^1$ piecewise open subset, then
$\displaystyle{\int_\gamma \omega=0}$. In the present example,
$U=\mathbb{R}^2\backslash\{(0,0)\}$ and let $\gamma$ be the unit
circle of parametric equations : $x_1(t)=\cos t$, $x_2(t)=\sin t$,
$t\in[0,2\pi]$. We have
$$
\int_\gamma\omega=\int_0^{2\pi}\left(-\frac{x_2}{x_1^2+x_2^2}x'_1(t)+\frac{x_1}{x_1^2+x_2^2}x'_2(t)\right)
=\int_0^{2\pi}(\sin^2t+\cos^2t)dt=2\pi.
$$
Since $\displaystyle{\int_\gamma\omega\neq0}$, then $\omega$ is
not exact This example shows that in Poincar\'{e} lemma, the
hypothesis that the open is starred is essential (here, the open
$U=\mathbb{R}^2\backslash\{(0,0)\}$ is not a star-shaped subset).
\end{rem}

\begin{prop}
The Poincar\'{e} lemma assures the existence of $\lambda$ but not
its uniqueness.
\end{prop}
\emph{Proof}: Indeed, let $\omega$ be a closed $k$-differential
form on a star-shaped open subset of $\mathbb{R}^n$. By
Poincar\'{e} lemma, there exist a $(k-1)$-differential form
$\lambda$ such that $\omega=d\lambda$. If $\mu$ is any
$(k-2)$-differential form, then $\lambda+d\mu$ satisfies the same
equation :
$$d(\lambda+d\mu)=d\lambda+d(d\mu)=d\lambda=\omega,$$
(because the exterior derivative obeys the rule : $d(d\mu)=0$, see
for example [7]). Conversely, if $\lambda_1$ and $\lambda_2$ are
any two $(k-2)$-differential forms such that :
$\omega=d\lambda_1=d\lambda_2$, then $d(\lambda_1-\lambda_2)=0$.
By Poincar\'{e} lemma, there exist a $(k-2)$-differential form
$\theta$ such that : $\lambda_1-\lambda_2=d\theta$, i.e.,
$\lambda_1=\lambda_2+d\theta$. From this we deduce that the
general solution can be expressed as the sum of a particular
solution and the derivative of an arbitrary $(k-2)$-differential
form. The proof is completed. $\square$

On manifods the Poincar\'{e} lemma can be stated as follows :

\begin{prop}
Any closed $k$-differential form $\omega$ is exact in the
neighborhood of an $n$-manifold $M$ (or, in $\mathbb{R}^n$ any
closed differential form is exact).
\end{prop}
\emph{Proof}: We have seen that for a star-shaped open subset of
$\mathbb{R}^n$, the form $\omega$ is exact. Since $M$ is locally
diffeomorphic to an open subset of $\mathbb{R}^n$, then for every
point $p\in M$ there exists also a neighborhood $U$ of $p$ and a
$(k-1)$-differential form $\lambda$ such that $\omega=d\lambda$ on
$U$. $\square$

\section{Some connections with Cech-De Rham-Dolbeault cohomologies}

We will give in the subsections below several formulations of
Poincar\'{e} lemma. Similarly, some examples and questions closely
related to the Poincar\'{e} lemma will be discussed.

\subsection{Cech-De Rham cohomologies and the Poincar\'{e} lemma}

Let $\mathcal{F}$ be a sheaf on a topological space $M$. Consider
the set $C^k(\mathcal{U},\mathcal{F})$ of $k$-cochains of degree
$k$ with values in $\mathcal{F}$, i.e., the set of applications
that associates with each $k$-up of open cover of $\mathcal{U}$ a
section on their intersection. In other words, we have
$$C^k(\mathcal{U},\mathcal{F})=
\prod_{\alpha_0,...,\alpha_k}\mathcal{F}(\mathcal{U}_{\alpha_0}\cap...\cap\mathcal{U}_{\alpha_k}),$$
where the direct product is taken over the set of all $(k+1)$
distinct elements $\alpha_0,...,\alpha_k$ in the index set. In
particular, we have
$C^0(\mathcal{U},\mathcal{F})=\displaystyle{\prod_\alpha\mathcal{F}(\mathcal{U}_\alpha)}$
and
$C^1(\mathcal{U},\mathcal{F})=\displaystyle{\prod_{\alpha,\beta}\mathcal{F}(\mathcal{U}_\alpha\cap\mathcal{U}_\beta)}$.
We define the coboundary operator
$$\delta:C^k(\mathcal{U},\mathcal{F})\longrightarrow
C^{k+1}(\mathcal{U},\mathcal{F}),$$ for $s\in
C^k(\mathcal{U},\mathcal{F})$, by
$$(\delta
s)(\mathcal{U}_{\alpha_0},...,\mathcal{U}_{\alpha_k})=
\sum_{j-0}^{k+1}(-1)^js\left.(\mathcal{U}_{\alpha_0},...,\widehat{\mathcal{U}_{\alpha_j}},...,
\mathcal{U}_{\alpha_{k+1}})\right|_{\mathcal{U}_{\alpha_0}\cap...\cap\mathcal{U}_{\alpha_{k+1}}}.$$
We easily check that $\delta^2=0$, so
$(C(\mathcal{U},\mathcal{F}))$ forms a cochain complex and we have
$$C^0(\mathcal{U},\mathcal{F})\overset{\delta }{\longrightarrow
}C^1(\mathcal{U},\mathcal{F})\overset{\delta }{\longrightarrow
}C^2(\mathcal{U},\mathcal{F})\overset{\delta }{\longrightarrow
}\cdots$$ with
\begin{eqnarray}
\delta&:&C^0(\mathcal{U},\mathcal{F})\longrightarrow
C^1(\mathcal{U},\mathcal{F}),\nonumber\\
&&s_{\alpha\beta}\longmapsto
s_\beta|_{\mathcal{U}_\alpha\cap\mathcal{U}_\beta}-s_\alpha|_{\mathcal{U}_\alpha\cap\mathcal{U}_\beta},\nonumber\\
\delta&:&C^1(\mathcal{U},\mathcal{F})\longrightarrow
C^2(\mathcal{U},\mathcal{F}),\nonumber\\
&&s_{\alpha\beta\gamma}\longmapsto
s_{\beta\gamma}|_{\mathcal{U}_\alpha\cap\mathcal{U}_\beta\cap\mathcal{U}_\gamma}
-s_{\alpha\gamma}|_{\mathcal{U}_\alpha\cap\mathcal{U}_\beta\cap\mathcal{U}_\gamma}
+s_{\alpha\beta}|_{\mathcal{U}_\alpha\cap\mathcal{U}_\beta\cap\mathcal{U}_\gamma},\nonumber\\
&\vdots&\nonumber\\
\delta&:&C^k(\mathcal{U},\mathcal{F})\longrightarrow
C^{k+1}(\mathcal{U},\mathcal{F}),\nonumber\\
&&s_{\alpha_0...\alpha_{k+1}}\longmapsto\sum_{j=0}^{k+1}(-1)^j\left.
s_{\alpha_0...\widehat{\alpha_j}...\alpha_{k+1}}\right|_{\mathcal{U}_{\alpha_0}\cap...\cap\mathcal{U}_{\alpha_{k+1}}}.\nonumber
\end{eqnarray}
For any sheaf $\mathcal{F}$ of abelian groups on $M$ and for any
open cover $\mathcal{U}$ of $M$, we can find an application
$\mathcal{F}\longrightarrow C^0(\mathcal{U},\mathcal{F})$, so that
the sequence
$$0\longrightarrow\mathcal{F}\longrightarrow C^0(\mathcal{U},\mathcal{F})\overset{\delta }{\longrightarrow
}C^1(\mathcal{U},\mathcal{F})\overset{\delta }{\longrightarrow
}\cdots$$ be exact. The Cech cohomology of the sheaf $\mathcal{F}$
with respect to the cover $\mathcal{U}$ is the quotient,
$$H^k(\mathcal{U},\mathcal{F})=\frac{\ker \left[\delta:C^k(\mathcal{U},\mathcal{F})\longrightarrow
C^{k+1}\mathcal{U},\mathcal{F})\right]}{\mbox{Im}
\left[\delta:C^{k-1}\mathcal{U},\mathcal{F})\longrightarrow
C^k\mathcal{U},\mathcal{F})\right]}.$$ These cohomology groups
depend on the cover $\mathcal{U}$. The Cech cohomology (or $k$-th
sheaf cohomology) of $\mathcal{F}$ on $M$ is defined to be the
direct limit of $H^k(\mathcal{U},\mathcal{F})$ as $U$ becomes
finer and finer,
$$H^k(M,\mathcal{F})=\lim_{\overrightarrow{\mathcal{U}}}H^k(\mathcal{U},\mathcal{F}).$$
The direct limit means, roughly, that for something to appear in
the final cohomology, it only needs to appear for sufficiently
refined covers. The passage to the limit in this definition is
delicate in practice although its interest is that there is no
longer any dependence on covers. However, Leray theorem [5]
asserts that :

\begin{thm}
If the open $\mathcal{U}_\alpha$ of a cover $\mathcal{U}$ of $M$
are such that :
$$H^k(A,\mathcal{F})=0,$$
for any $k>0$ and any finite intersection
$A\equiv\mathcal{U}_{\alpha_1}\cap...\cap\mathcal{U}_{\alpha_k}$
of open $\mathcal{U}_\alpha$, then the groups
$H^k(\mathcal{U},\mathcal{F})$ et $H^k(M,\mathcal{F})$ are
isomorphic for all $k$.
\end{thm}

Let $\Omega^k$ be the sheaf of $k$-differential forms on a
manifold $M$ and $\mathcal{Z}^k$ be the sheaf of the closed
$k$-forms. Let $d_*$ be the application of cohomology associated
with $d : \Omega^k\longrightarrow\mathcal{Z}^{k+1}$. The $k$-th De
Rham cohomology $H^k_{DR}(M)$ is defined as the quotient space of
the $k$-closed differential forms by the $(k-1)$-differential
forms. Since a $k$-form is a global section of the sheaf of
$k$-forms, so it is an element of $H^0(M, \Omega^k)$ because this
latter group of cohomology identifies with global sections.
Therefore,
$$H^k_{DR}(M)=\frac{H^0(M, \mathcal{Z}^k)}{d_*H^0(M,
\Omega^{k-1})}.$$ As a consequence of the Poincar\'{e} lemma, we
have the following result :

\begin{prop}
De Rham's cohomology of a manifold is isomorphic to its Cech
cohomology with coefficients in $\mathbb{R}$.
\end{prop}
\emph{Proof}: According to Poincar\'{e} lemma, any closed form is
locally exact. So the sheaf sequence
$$0\longrightarrow\mathcal{Z}^k\longrightarrow \Omega^k
\overset{d}{\longrightarrow }\mathcal{Z}^{k+1}\longrightarrow 0,$$
is exact where $\mathcal{Z}^0\equiv \mathbb{R}$ is the sheaf of
locally constant functions and $\Omega^0$ the sheaf of fonctions
in $\mathcal{C}^\infty$. The following long exact sequence in
cohomology associated with the above sequence is
$$H^q(M, \Omega^k)\overset{d^*}{\longrightarrow} H^q(M,
\mathcal{Z}^{k+1})\overset{\partial}{\longrightarrow} H^{q+1}(M,
\mathcal{Z}^k)\longrightarrow H^{q+1}(M, \Omega^k),$$ where $d^*$
denotes the exterior derivative and $\partial$ the boundary
operator of the long exact sequence. The sheaf $\Omega^p$ admits
partitions of the unit, hence for $q\geq 1$, $$H^q(M,
\Omega^k)=0,$$ and so
$$H^q(M, \mathcal{Z}^{k+1})=H^{q+1}(M, \mathcal{Z}^k).$$
For the particular case $q=0$, we have
$$H^0(M, \Omega^k)\overset{d^*}{\longrightarrow} H^0(M,
\mathcal{Z}^{k+1})\overset{\partial}{\longrightarrow} H^1(M,
\mathcal{Z}^p)\longrightarrow 0,$$ and
$$H^1(M, \Omega^k)=\frac{H^0(M, \mathcal{Z}^{k+1})}{d_*H^0(M,
\Omega^k)}.$$ Therefore,
$$
H^k_{DR}(M)=H^1(M, \mathcal{Z}^{k-1})=H^{k-1}(M,
\mathcal{Z}^1)=H^k(M, \mathcal{Z}^0)=H^k(M, \mathbb{R}),$$ and the
result follows. $\square$

Let $M$ and $N$ two manifolds. Recall that two maps
$f:M\longrightarrow N$ and $g:N\longrightarrow M$ are (smoothly)
homotopic, if there exists a (smooth) homotopy
$$h:M\times
I\longrightarrow N,\qquad h(x,0)=f(x),\qquad h(x,1)=g(x),$$ ($I$
is the interval, $0\leq t\leq1$). We say the manifolds $M$ and $N$
are homotopy equivalent if there exist (smooth) maps
$f:M\longrightarrow N$ and $g:N\longrightarrow M$ such that the
composites $gof:M\longrightarrow M$ and $fog:N\longrightarrow N$
are (smoothly) homotopic to the respective identity maps
$M\longrightarrow M$, $x\longmapsto x$ and $N\longrightarrow N$,
$y\longmapsto y$. For example, the closed unit disc in
$\mathbb{R}^n$ is homotopy equivalent to a point. It easy to
verify that $\mathbb{R}^n\backslash\{0\}$ and $S^{n-1}$ are
homotopy equivalent. Using reasoning similar to proof 2, we show
that :

\begin{prop}
If two manifolds $M$ and $N$ are homotopy equivalent then their
cohomology groups are isomorphic.
\end{prop}
\emph{Proof}: Indeed, consider the (smooth) homotopy $F:M\times
I\longrightarrow N$ between the above two maps $f$ and $g$. Then
the induced homomorphisms
\begin{eqnarray}
f^*&:&H^k(N)\longrightarrow
H^k(M),,\nonumber\\
g^*&:&H^k(M)\longrightarrow H^k(N),\nonumber
\end{eqnarray}
of the cohomology groups coincide. It suffices to consider a
differential form $\omega\in H^k(gof)(M)$, which implies that it
exists a pullback
\begin{eqnarray}
h^*(\omega)&=&\lambda+dt\wedge\theta,\nonumber\\
h^*(\omega)|_{t=t_0}&=&h^*(.,t_0)(\omega)=\lambda|_{t=t_0},\nonumber
\end{eqnarray}
where $\lambda\in H^k(M\times I)$, $\theta\in H^{k-1}(M\times I)$,
and which $\lambda$, $\theta$ do not involve the differential $dt$
(in the sens that $\lambda$, $\theta$ do not contain a $dt$ term).
We have\footnote{Let $V\subset \mathbb{R}^m$, $U\subset
\mathbb{R}^n$ be open subsets, $g : V\longrightarrow U$ be a
differentiable application and $\displaystyle{\omega=\sum_{1\leq
i_1<...<i_k\leq n}f_{i_1,...,i_k}dx_{i_1}\wedge...\wedge
dx_{i_k}}$, a $k$-differential form on $U$. If $\omega$ is
$\mathcal{C}^1$ on $U$ and $g$ is $\mathcal{C}^2$ on $V$, then
(see for example [8]), $g^*(d\omega)=d(g^*\omega)$.}
$$h^*(d_M(\omega))=d_{M\times I}(h^*(\omega))=dt\wedge\left(\frac{\partial\lambda}{\partial
t}-d_M(\theta)\right),$$ and since $\omega$ is closed then
$h^*(d_M(\omega))=0$, so
$$\frac{\partial\lambda}{\partial t}=d_M(\theta).$$
Moreover, we have
\begin{eqnarray}
h^*(.,1)(\omega)-h^*(.,0)(\omega)&=&\lambda|_{t=1}-\lambda|_{t=0},\nonumber\\
&=&\int_0^1\frac{\partial\lambda}{\partial
t}dt,\nonumber\\
&=&\int_0^1d_M(\theta)dt,\nonumber\\
&=&d_M\left(\int_0^1\theta dt\right).\nonumber
\end{eqnarray}
Therefore,
$$(gof)^*=h^*(.,1),$$
is the identity map hence, $f^*$ and $g^*$ are inverse to each
other and in particular, both are isomorphisms. More precisely,
since $gof\sim id_M$ and $fog\sim id_N$, we deduce that
\begin{eqnarray}
f^*og^*&=&(gof)^*=id^*_M=id_{H^k(M)},\nonumber\\
g^*of^*&=&(fog)^*=id^*_N=id_{H^k(N)}.\nonumber
\end{eqnarray}
Hence, $f^*$ is a vector space isomorphism and $(f^*)^{-1}=g^*$,
which completes the proof. $\square$

\begin{rem}
Recall that a (smooth) manifold $M$ is contractible if the
identity map on $M$ is homotopic to a constant map or in other
words if $M$ is homotopy equivalent to a point. For example, a
star-shaped open subset of $\mathbb{R}^n$ is contractible and
hence has De Rham cohomology of a point. The cohomology groups of
space $\mathbb{R}^n$ (and of the ball around any point $p$) are
isomorphic to those of $p$. Thus $H^k(\mathbb{R}^n)$ is trivial
for $k>0$, while $H^0(\mathbb{R}^n)\simeq\mathbb{R}$. In fact if
$M$ is contractible, which means that $M$ is homotopy equivalent
to a point and since homotopy equivalent manifolds admit
isomorphic De Rham cohomology groups, then
$$H^k(M)=\left\{\begin{array}{rl}
\mathbb{R},& k=0\\
0,&  k>0
\end{array}\right.$$
This fact leads immediately to the Poincar\'{e} Lemma. As we
mentioned above when $k=0$, notice that constant functions are
closed and $H^0(M)$ is a finite dimensional vector space equal to
the number of connected components of $M$. Here we have only
$0$-differential forms, i.e., $ f(x)$ functions on $M$. There are
no exact differential forms. Then, $H^0(M,\mathbb{R})=\{f: f
\mbox{ is closed}\}$. Since $df(x)=0$, then in any chart
$(U,x_1,...,x_n)$ of $M$, we have
$$\frac{\partial f}{\partial x_1}=\frac{\partial f}{\partial
x_2}=\cdots=\frac{\partial f}{\partial x_n}=0.$$ Therefore,
$f(x)=$ constant locally, i.e., $f(x)=$ constant on each connected
component of $M$. So the number of connected components of  $M$ is
the dimension in question.
\end{rem}

\begin{rem}
As we have seen the Poincar\'{e} lemma and De Rham theorem give a
characterization of the de Rham cohomology groups. Applying
Poincar\'{e} lemma for $k$-differential forms where $k\neq0$, we
show that the cohomology of a ball is trivial and in particular,
$\dim H^k(\mathbb{R}^n)=\delta(k, 0)$. If $S^n$ is the $n$-sphere
and if $T^n=(S^1)^n$ is the $n$-dimensional torus, i.e., the
product of $n$ circles, then by proceeding by induction on $n$, we
obtain $$\dim H^k(S^n)=\delta(n, k),\qquad \dim
H^k(T^n)=\frac{n!}{k!(n-k)!}.$$
\end{rem}

\subsection{Dolbeault cohomology : $\overline{\partial}$-Poincar\'{e} lemma (or
Dolbeault-Grothendieck lemma)}

In each point of a complex manifold $M$, we define local
coordinates $z_j$ and $\overline{z}_j$ as well as the tangent
bundle $\mathbb{C}\left[\frac{\partial}{\partial z_j},
\frac{\partial}{\partial \overline{z}_j}\right]$ that we note
$TM$. This admits the decomposition $TM=T'M\oplus T''M$, where
$T'M$ (resp. $T''M$) is the holomorphic (respectively
antiholomorphic) part of $TM$ generated by
$\displaystyle{\frac{\partial}{\partial z_j}}$ (resp.
$\displaystyle{\frac{\partial}{\partial \overline{z}_j}}$).
Similarly, the cotangent bundle $T^*M$ (dual of $TM$) admits the
decomposition : $T^*M=T'^*M\oplus T''^*M$, into holomorphic and
antiholomorphic parts. The points of this bundle are the linear
forms on the fibers of $TM\longrightarrow M$. A $p$-differential
form is a section of the bundle $\Lambda^pT^*M$ (the exterior
powers of $T^*M$). The points of the latter being the alternating
multilinear $p$-forms on $T^*M$. If $z_1,...,z_n$ are holomorphic
local coordinates on a complex manifold $M$ of dimension $n$, then
$dz_1,...,dz_n, d\overline{z}_1,...,d\overline{z}_n$ form a local
base of the space of differential forms. A differential form on
$M$ of type $(p, q)$ is given locally by
$$\omega=\sum_{\overset{1\leq j_1<...<j_p\leq n}{1\leq k_1<...<k_q\leq
n}}f_{j_1...j_pk_1...k_q}dz_{j_1}\wedge...\wedge dz_{j_p}\wedge
d\overline{z}_{k_1}\wedge...\wedge d\overline{z}_{k_q},$$ where
$f_{j_1...j_pk_1...k_q}$ are $\mathcal{C}^\infty$ functions with
complex values. This definition does not depend on the choice of
holomorphic local coordinates. Obviously, a form of type $(0,0)$
is a function. The set of differential forms of type $(p, q)$ over
$M$ is a complex vector bundle denoted $A^{p, q}_M$. It should be
noted that only $A^{p, 0}_M$ is holomorphic. The differential of
$\omega$ is
\begin{eqnarray}
d\omega&=&\sum_{\overset{1\leq j_1<...<j_p\leq n}{1\leq
k_1<...<k_q\leq n}}df_{j_1...j_pk_1...k_q}\wedge
dz_{j_1}\wedge...\wedge dz_{j_p}\wedge
d\overline{z}_{k_1}\wedge...\wedge
d\overline{z}_{k_q},\nonumber\\
&=&\sum_{\overset{1\leq j_1<...<j_p\leq n}{1\leq k_1<...<k_q\leq
n}}\partial f_{j_1...j_pk_1...k_q}\wedge dz_{j_1}\wedge...\wedge
dz_{j_p}\wedge d\overline{z}_{k_1}\wedge...\wedge
d\overline{z}_{k_q}\nonumber\\
&&\quad +\sum_{\overset{1\leq j_1<...<j_p\leq n}{1\leq
k_1<...<k_q\leq n}}\overline{\partial}f_{j_1...j_pk_1...k_q}\wedge
dz_{j_1}\wedge...\wedge dz_{j_p}\wedge
d\overline{z}_{k_1}\wedge...\wedge
d\overline{z}_{k_q},\nonumber\\
&=&\partial \omega+\overline{\partial}\omega,\nonumber
\end{eqnarray}
where $\displaystyle{\partial\equiv\frac{\partial}{\partial z}dz}$
and
$\displaystyle{\overline{\partial}\equiv\frac{\partial}{\partial
\overline{z}}d\overline{z}}$. Since $\partial
f_{j_1...j_pk_1...k_q}$ is of type $(1,0)$ and
$\overline{\partial} f_{j_1...j_pk_1...k_q}$ is of type $(0,1)$,
then the component $\partial \omega$ of $d\omega$ is of type
$(p+1,q)$ and $\overline{\partial}\omega$ is of type $(p,q+1)$. So
we have the following decomposition :
$d=\partial+\overline{\partial}$ and the operators $\partial$,
$\overline{\partial}$ are independent of the choice of holomorphic
local coordinates. The relation $d^2=0$ means here that :
$$\partial^2=\overline{\partial}^2=\partial\overline{\partial}+\overline{\partial}\partial=0.$$
Indeed, we deduce from the relations $d^2=0$ and
$d=\partial+\overline{\partial}$, that
$$0=d^2=\partial^2+\overline{\partial}^2+\partial\overline{\partial}+\overline{\partial}\partial.$$
Note that if $\omega$ is a form of type $(p, q)$, then
$\partial^2\omega$ is of type $(p+2,q)$,
$\overline{\partial}^2\omega$ is of type $(p,q+2)$ and
$(\partial\overline{\partial}+\overline{\partial}\partial)$ is of
type $(p+1,q+1)$. The result is that the sum of these three forms
is zero. The operator $\overline{\partial}$ satisfies the
following rule (Leibniz rule) :
$$\overline{\partial}(\omega\wedge\lambda)=\overline{\partial}\omega\wedge\lambda+
(-1)^{{deg}\omega}\omega\wedge\overline{\partial}\lambda.$$
Noticing that
$\partial\omega=\overline{\overline{\partial}\overline{\omega}}$,
we also obtain
$$\partial(\omega\wedge\lambda)=\partial\omega\wedge\lambda+
(-1)^{{deg }\omega}\omega\wedge\partial\lambda.$$ Moreover, the
differential form $\omega$ is said $\overline{\partial}$-closed if
$\overline{\partial}\omega=0$ and $\overline{\partial}$-exact if
there exists a form $\lambda$ such that :
$\omega=\overline{\partial}\lambda$. Any
$\overline{\partial}$-exact form is $\overline{\partial}$-closed.
The analogue of Poincar\'{e} lemma for the Dolbeault operator
$\overline{\partial}$ (see below for definition) is the
$\overline{\partial}$-Poincar\'{e} lemma (or
Dolbeault-Grothendieck lemma) and can be stated as follows :
\begin{prop}
If $\omega$ is a $\mathcal{C}^1$ form of type $(p,q)$, $q>0$ with
$\overline{\partial}\omega=0$, then it exists locally on $M$ a
$\mathcal{C}^1$ form $\lambda$ of type $(p,q-1)$ such that :
$\omega=\overline{\partial}\lambda$.
\end{prop}

We have shown previously that Poincar\'{e} lemma leads to the
nullity of cohomology groups in the case of star-shaped open
subsets. We have an isomorphism
$$H^{p,q}_{\overline{\partial}}(M)\simeq H^q(M, \Omega_M^p),$$ for
all integers $p$ and $q$, where
$$H^{p,q}_{\overline{\partial}}(M)=\frac{\{\overline{\partial}\mbox{-closed forms of type }(p,q)\mbox{ on } M\}}
{\{\overline{\partial}\mbox{-exact forms of type }(p,q) \mbox{ on
} M\}},$$ is a Dolbeault cohomology group of complex manifold $M$.


\begin{thebibliography}{99}

\bibitem{}V. I. Arnold, \emph{Mathematical methods in classical
mechanics}, Springer-Verlag, Berlin-Heidelberg- New York, 1978.
\bibitem{}V. I. Arnold, A.B. Givental, \emph{Symplectic geometry in :
Dynamical systems IV}, (eds). V.I. Arnold and S.P. Novikov, (EMS,
Volume \textbf{4}, pp. 1-136) Springer-Verlag, 1988.
\bibitem{} H. Cartan, \emph{Formes différentielles. Application élémentaires au calcul des variations at à la théorie des courbes et surfaces},
Hermann, Paris, 2007.
\bibitem{}B. A. Dubrovin, S.P., Novikov, A.T., Fomenko, A.T., \emph{Modern geometry.
Methods and applications}, Parts I, II, Springer-Verlag, 1984,
1985.
\bibitem{}P.A. Griffiths, J. Harris, \emph{Principles of algebraic
geometry}, Wiley-Interscience, New-York, 1978.
\bibitem{} S. Hans, \emph{Differential Forms, the Early Days; or the Stories
of Deahna's Theorem and of Volterra's Theorem}, The American
Mathematical Monthly, Vol. 108, No. 6 (2001), 522-530
\bibitem{} A. Lesfari, \emph{Introduction \`{a} la g\'{e}om\'{e}trie alg\'{e}brique complexe}, \'{E}ditions Hermann, Paris,
2015.
\bibitem{} A. Lesfari, \emph{Formes différentielles et analyse vectorielle (Cours et exercices résolus)},
éditions Ellipses, Paris, 2017.
\bibitem{} J. Mawhin, \emph{Analyse. Fondements, techniques, évolution}, De Boeck Université, Bruxelles, 1997.
\bibitem{} V. Volterra, \emph{Delle variabili complesse negli iperspazii}, Rend.
Accad. dei Lincei, ser. IV, vol. V1 (1889), NotaI, pp. 158-165,
Nota II, pp.291-299; \emph{Sulle funzione conjugate}, ibidem, ser.
Iva, vol. VI (1889), p. 158-169; \emph{Opere Matematiche}, Accad.
Nazion. dei Lincei, Roma (1954), vol. 1 , pp. 403-410, 411-419,
420-432.
\bibitem{}A. Weinstein, \emph{Symplectic manifolds and their lagrangian submanifolds},
Advances in Maths., 6, (1971), 329-346.
\bibitem{}A. Weinstein, \emph{Lectures on ymplectic manifolds},
American Mathematical Society, Cbms Regional Conference Series in
Mathematics, Number 29, 1977.
\bibitem{}A. Weinstein, \emph{The local structure of Poisson manifolds}, J. Diff.
Geom., 18, (1983), 523-557.

\end{thebibliography}
\end{document}